\documentclass[reqno]{amsart}
\usepackage[utf8]{inputenc}
\usepackage{amsmath}
\usepackage{amssymb}
\usepackage{amsthm}
\usepackage[all]{xy}
\xyoption{tips}
\SelectTips{cm}{10}
\usepackage[dvipsnames,table,xcdraw]{xcolor}
\usepackage{fullpage}
\usepackage{stmaryrd}
\usepackage{enumitem}
\usepackage{mathtools}
\usepackage{microtype}
\usepackage{todonotes,booktabs}
\usepackage{float}
\usepackage{aliascnt}

\usepackage{tabularray}

\usepackage[parfill]{parskip}

  \newcommand{\Ho}{\mathrm{Ho}}
  
\usepackage{tikz,scalerel}
\usetikzlibrary{positioning,arrows}
\usetikzlibrary{shapes.geometric, calc}
\usetikzlibrary{decorations.pathreplacing}

\theoremstyle{plain}

\newaliascnt{theorem}{equation}  
  
\aliascntresetthe{theorem}

\theoremstyle{definition}

\newaliascnt{proposition}{equation}  

\aliascntresetthe{proposition}

 \newaliascnt{hproposition}{equation}  

\aliascntresetthe{hproposition}

\newaliascnt{lemma}{equation}  

\aliascntresetthe{lemma}

\newaliascnt{corollary}{equation}  

\aliascntresetthe{corollary}

\newaliascnt{claim}{equation}  

\aliascntresetthe{claim}

\newaliascnt{conjecture}{equation}  

\aliascntresetthe{conjecture}

\newaliascnt{question}{equation}  

\aliascntresetthe{question}

\newaliascnt{definition}{equation}  
\newtheorem{definition}[definition]{Definition}
\aliascntresetthe{definition}

\newaliascnt{hdefinition}{equation}  

\aliascntresetthe{hdefinition}

\newaliascnt{example}{equation}  

\aliascntresetthe{example}

\theoremstyle{remark}

\newaliascnt{remark}{equation}  

\aliascntresetthe{remark}

\newaliascnt{convention}{equation}  

\aliascntresetthe{convention}

\usepackage{hyperref}

\definecolor{dark-red}{rgb}{0.5,0.15,0.15}
\definecolor{dark-blue}{rgb}{0.15,0.15,0.6}
\definecolor{dark-green}{rgb}{0.15,0.6,0.15}
\hypersetup{
    colorlinks, linkcolor=codepurple,
    citecolor=dark-blue, urlcolor=magenta
}

\newcommand{\sfT}{\mathsf{T}}

\usepackage{mathrsfs}

\usepackage{bbm}
\usepackage{multicol}

\usepackage{listings}
\usepackage{xcolor}

\definecolor{codegreen}{rgb}{0,0.6,0}
\definecolor{codegray}{rgb}{0.5,0.5,0.5}
\definecolor{codepurple}{rgb}{0.58,0,0.82}
\definecolor{backcolour}{rgb}{0.95,0.95,0.92}

\lstdefinestyle{mystyle}{
    backgroundcolor=\color{backcolour},   
    commentstyle=\color{codegreen},
    keywordstyle=\color{magenta},
    numberstyle=\tiny\color{codegray},
    stringstyle=\color{codepurple},
    basicstyle=\ttfamily\footnotesize,
    breakatwhitespace=false,         
    breaklines=true,                 
    captionpos=b,                    
    keepspaces=true,                 
    numbers=left,                    
    numbersep=5pt,                  
    showspaces=false,                
    showstringspaces=false,
    showtabs=false,                  
    tabsize=2
}

\makeatletter
\def\@tocline#1#2#3#4#5#6#7{\relax
  \ifnum #1>\c@tocdepth % then omit
  \else
    \par \addpenalty\@secpenalty\addvspace{#2}%
    \begingroup \hyphenpenalty\@M
    \@ifempty{#4}{%
      \@tempdima\csname r@tocindent\number#1\endcsname\relax
    }{%
      \@tempdima#4\relax
    }%
    \parindent\z@ \leftskip#3\relax \advance\leftskip\@tempdima\relax
    \rightskip\@pnumwidth plus4em \parfillskip-\@pnumwidth
    #5\leavevmode\hskip-\@tempdima
      \ifcase #1
       \or\or \hskip 1em \or \hskip 2em \else \hskip 3em \fi%
      #6\nobreak\relax
    \dotfill\hbox to\@pnumwidth{\@tocpagenum{#7}}\par
    \nobreak
    \endgroup
  \fi}
\makeatother

\usepackage[nameinlink,capitalise,noabbrev]{cleveref}

\title{ninfty: A software package for homotopical combinatorics}

 \author{Scott Balchin}
 \address{Mathematical Sciences Research Centre, Queen's University Belfast, UK}
  \email{s.balchin@qub.ac.uk }

\begin{document}

\begin{abstract}
We introduce \lstinline{ninfty}\thanks{\url{https://github.com/bifibrant/ninfty}}, a header-only \lstinline{C++} library distributed under an MIT Open Source License designed for the study of enumeration problems arising in homotopical combinatorics. The \lstinline{ninfty} repository moreover contains a folder with data files for many common finite groups. This is in addition to Sage code which can be used to generate input data for further finite groups, and Sage code for generating input data for abstract lattices which may not arise at the subgroup lattice of a group.
\end{abstract}

\maketitle

\contentsline {section}{\tocsection {}{1}{Context}}{2}{section.1}%
\contentsline {section}{\tocsection {}{2}{Installation, requirements, and limitations}}{3}{section.2}%
\contentsline {section}{\tocsection {}{3}{Minimal working example of \lstinline {ninfty}}}{3}{section.3}%
\contentsline {section}{\tocsection {}{4}{Generating data files using Sage}}{5}{section.4}%
\contentsline {section}{\tocsection {}{5}{Feature description}}{6}{section.5}%
\contentsline {subsection}{\tocsubsection {}{5.1}{Basic features}}{6}{subsection.5.1}%
\begin{multicols}{2}
\parskip=0pt \relax
\contentsline {subsubsection}{\lstinline {printNumberOfTransfers}}{6}{lstnumber.-10.1}%
\contentsline {subsubsection}{\lstinline {printNumberOfCosaturatedTransfers}}{6}{lstnumber.-11.1}%
\contentsline {subsubsection}{\lstinline {printNumberOfSaturatedTransfers}}{6}{lstnumber.-12.1}%
\contentsline {subsubsection}{\lstinline {printNumberOfUnderlyingTransfers}}{7}{lstnumber.-13.1}%
\contentsline {subsubsection}{\lstinline {printNumberOfConjugacyTransfers}}{7}{lstnumber.-14.1}%
\contentsline {subsubsection}{\lstinline {printWidth}}{7}{lstnumber.-15.1}%
\contentsline {subsubsection}{\lstinline {printComplexity}}{8}{lstnumber.-16.1}%
\contentsline {subsubsection}{\lstinline {printNumberOfMaximallyGenerated}}{8}{lstnumber.-17.1}%
\contentsline {subsubsection}{\lstinline {printNumberOfFlatTransfers}}{8}{lstnumber.-18.1}%
\contentsline {subsubsection}{\lstinline {printNumberOfTransferPairs}}{9}{lstnumber.-19.1}%
\contentsline {subsubsection}{\lstinline {printNumberOfCClosedPairs}}{9}{lstnumber.-20.1}%
\contentsline {subsubsection}{\lstinline {printNumberOfQuillenPairs}}{9}{lstnumber.-21.1}%
\contentsline {subsubsection}{\lstinline {printNumberOfWeakEquivalenceTypes}}{10}{lstnumber.-22.1}%
\contentsline {subsubsection}{\lstinline {printNumberOfCompatiblePairs}}{10}{lstnumber.-23.1}%
\contentsline {subsubsection}{\lstinline {printSageTransferPoset}}{10}{lstnumber.-24.1}%
\contentsline {subsubsection}{\lstinline {printSageCClosedPoset}}{10}{lstnumber.-25.1}%
\contentsline {subsubsection}{\lstinline {printSageQuillenPoset}}{10}{lstnumber.-26.1}%
\contentsline {subsubsection}{\lstinline {dataSheet} and \lstinline {dataSheetLatex}}{11}{lstnumber.-27.2}%
\contentsline {subsubsection}{\lstinline {printSubgroupDictionary}}{11}{lstnumber.-28.1}%
\contentsline {subsubsection}{\lstinline {printAllTransfers}}{12}{lstnumber.-30.1}%
\vfill\null
\end{multicols}
\contentsline {subsection}{\tocsubsection {}{5.2}{Advanced features}}{12}{subsection.5.2}%
\begin{multicols}{2}
\parskip=0pt \relax
\contentsline {subsubsection}{\lstinline {transferFind}}{12}{lstnumber.-32.1}%
\contentsline {subsubsection}{\lstinline {isSaturated}}{13}{lstnumber.-33.1}%
\contentsline {subsubsection}{\lstinline {isCosaturated}}{13}{lstnumber.-34.1}%
\contentsline {subsubsection}{\lstinline {isFlat}}{13}{lstnumber.-35.1}%
\contentsline {subsubsection}{\lstinline {leftSet}}{13}{lstnumber.-36.1}%
\contentsline {subsubsection}{\lstinline {saturatedHull}}{13}{lstnumber.-37.1}%
\contentsline {subsubsection}{\lstinline {cosaturatedCore}}{13}{lstnumber.-38.1}%
\contentsline {subsubsection}{\lstinline {dualTransferSystem}}{13}{lstnumber.-39.1}%
\contentsline {subsubsection}{\lstinline {minimalFibrantSubgroup}}{14}{lstnumber.-40.1}%
\contentsline {subsubsection}{\lstinline {findBasis}}{14}{lstnumber.-41.1}%
\contentsline {subsubsection}{\lstinline {weakEquivalences}}{14}{lstnumber.-42.1}%
\contentsline {subsubsection}{\lstinline {modelCheck}}{14}{lstnumber.-43.1}%
\contentsline {subsubsection}{\lstinline {isCompatible}}{14}{lstnumber.-44.1}%
\contentsline {subsubsection}{\lstinline {edgesToTikz}}{14}{lstnumber.-45.1}%
\vfill\null
\end{multicols}
\contentsline {section}{\tocsection {}{6}{Support}}{16}{section.6}%
\contentsline {section}{\tocsection {}{}{References}}{17}{section*.2}%

\newpage

\section{Context}\label{sec:context}

\textit{Homotopical combinatorics} is a rapidly developing field that is concerned with combinatorial phenomena that arises in the study of a fruitful intersection of equivariant homotopy theory, equivariant algebra, and abstract homotopy theory.

The theory begins its journey in the realm of studying equivariant gadgets equipped with a homotopically meaningful multiplicative structure. Non-equivariantly, homotopy commutativity is governed by $E_\infty$-operads, and up to homotopy, we are  agnostic to the choice of a $E_\infty$-operad, whether that be the Barrett--Eccles operad or the little $\infty$-disks operad. That is, all $E_\infty$-operads have equivalent categories of algebras \cite{iterated}.

In the presence of the action of a finite group $G$ the situation is more exotic. In essence, one wants to encode how much of the group action is respected by the commutativity of a multiplication. At one extreme we may consider the group virtually redundant, and we are proffered \emph{naive} commutative objects. At the other extreme, we find ourselves with \emph{genuine} commutative objects. In particular, we have non-equivalent $G$-equivariant generalizations of $E_\infty$-operads. This observation led Blumberg--Hill to define the notion of \emph{$N_\infty$-operads}~\cite{blumberghill}, whence the namesake of this computational package.

The wealth of interest in this theory comes from the quagmire that occupies the space in-between the aforementioned extremes. We are primary concerned in studying the finite poset $\Ho(N_\infty(G))$ of $N_\infty$-operads for a fixed group $G$ up to homotopy. This gives us an insight into the multitudes of ways in which a given object may be homotopy commutative with respect to a group action. However, the definition of an $N_\infty$-operad is somewhat complex, and not amenable to investigative study.

A key step in realizing this goal comes from work of Rubin \cite{Rubin_comb} and Balchin--Barnes--Roitzheim \cite{bbr}, building on work of Blumberg--Hill \cite{blumberghill}, Bonventre--Pereira \cite{BP_operads}, and Guti\'errez--White \cite{GW_operads}, which extracts a discrete combinatorial gadget called a \emph{$G$-transfer system} such that the poset of $G$-transfer systems, $\mathsf{Tr}(G)$, is order isomorphic to the poset $\Ho(N_\infty(G))$. Explicitly:

\begin{definition}\label{defn:xfer}
    Let $G$ be a finite group. A \textit{$G$-transfer system} $\mathsf{T}$ is a collection of pairs of subgroups $(H,K)$ with $1 \leqslant H \leqslant K \leqslant G$ such that
    \begin{itemize}
        \item (Identity) $(H,H) \in \mathsf{T}$ for all $H \leqslant G$;
        \item (Composition) If $(L,K) \in \mathsf{T}$ and $(K,H) \in T$ then $(L,H) \in \mathsf{T}$;
        \item (Conjugation) If $(K,H) \in \mathsf{T}$ then $(gKg^{-1} , gHg^{-1}) \in \mathsf{T}$ for all $g \in G$;
        \item (Restriction) If $(K,H) \in \mathsf{T}$ and $L \leqslant H$ then $(K \cap L, L) \in \mathsf{T}$.
    \end{itemize}
\end{definition}
As such, a transfer system can be described as a subgraph of the lattice of subgroups with inclusion $\operatorname{Sub}(G)$, with $G$ acting on this lattice via conjugacy.

The definition of a transfer system is combinatorially simple, and thus we have a chance to explore the collection of them for a finite group $G$. This, however, is just the tip of the iceberg of homotopical combinatorics. Work of many authors has developed an extensive theory which yields surprising links between seemingly disparate areas of mathematics. We refer the reader to \cite{zbMATH07918757} for a comprehensive overview of the current landscape of homotopical combinatorics.

Many of the developments in this area have been informed initially by computational evidence. The goal of \lstinline{ninfty} is to provide a collection of high performance implementations of all algorithms required to generate this computational data for known constructions of interest, removing a potential barrier to entry in this rapidly developing field.

This document has two primary purposes. First, it serves as a user guide to \lstinline{ninfty}, and provides a minimal working example of the complete pipeline of obtaining results in \cref{sec:MWE} and \cref{sec:Sage}. Secondly, it lists the features contained in \lstinline{ninfty}, along with the references of where the concept was introduced, or has been studied, in the literature, along with an exemplar result. It has been written in such a way that even those who are not confident in using \lstinline{C++} can avail of the features.

\newpage

The description of the functionality has been split in two distinct sections:
\begin{enumerate}[align=left]
    \item[(Basic)] Described in \cref{sec:featuresBASIC}, these are features which print a requested value for a given group to the console and require no knowledge of \lstinline{C++}.
    \item[(Advanced)] Described in \cref{sec:featuresADVANCED}, these features allow one to explicitly access the stored variables and look at individual objects. Along with these advanced features are functions which can be used to produce TikZ representations of various structures. The use of these features requires a base knowledge of \lstinline{C++}, of which many tutorials exist online.
\end{enumerate}

Many of the features presented also make sense for lattices thanks to work of \cite{fooqw}, which proves for Dedekind groups that transfer systems are in bijection with weak factorization systems on the underlying lattice of subgroups. As such, \lstinline{ninfy} also allows for the input to be an arbitrary lattice.

\textbf{Acknowledgments.} The author would like to thank Kyle Ormsby and Ben Spitz for many enlightening discussions about the development of \lstinline{ninfty}, along with many helpful suggestions for improvements of this document. Moreover, the author would like to thank the participants of the AMS MRC on Homotopical Combinatorics and the participants of the eCHT Reading Group on Homotopical Combinatorics for their enthusiasm and willingness to test some features of \lstinline{ninfty}.

\section{Installation, requirements, and limitations}\label{sec:install}

To run \lstinline{ninfty}, only \lstinline{C++20} is required, no further libraries are needed. The code that is used for generation of data for a given group or lattice requires an installation of Sage.

The first step is to clone or download the \lstinline{ninfty} repository from
\url{https://github.com/bifibrant/ninfty}. There are four files of interest in this repository:
\begin{itemize}
    \item \lstinline{ninfty.h}: This is the \lstinline{C++} header file which contains all the functionality of \lstinline{ninfty}.
    \item \lstinline{main.cpp}: This is an example driver file. Note that it \lstinline{includes} the \lstinline{ninfty} header file along with a \lstinline{group_data/G.h}, where $G$ is some finite group. The folder \lstinline{group_data} contains data files for common finite groups of interest. It is in this driver file that functions should be run (see \cref{sec:MWE} for a minimal working example).
    \item \lstinline{group_gen.sage}: A Sage file, located in the \lstinline{group_data} folder which can be used to generate further header files for finite groups.
    \item \lstinline{lattice_gen.sage}: A Sage file, located in the \lstinline{lattice_data} folder which can be used to generate further header files for finite lattices.
\end{itemize}

Many of the answers to the enumeration problems seem to scale quickly with the number of (conjugacy classes of) subgroups. Given the finite limitations of computing, this means that it is not feasible to expect to be able to compute these invariants for a large class of finite groups, and some algorithms scale worse than others (for example, algorithms which ask about certain pairs of transfer systems). As a benchmark, computing all 5,389,480 transfer systems for $C_{pqrs}$ takes approximately 450 seconds on a base 2023 Macbook Pro. 

The main limitation with \lstinline{ninfty} is with memory, in particular the algorithm to generate all transfer systems must store all of them as it runs. While the data is stored in a lightweight form, the scale of the computations can lead to difficulty. For example the number of transfer systems for $S_5$ (183,598,202) was achieved only using HPC facilities as it required upwards of 100GB of RAM to complete.

\section{Minimal working example of \lstinline{ninfty}}\label{sec:MWE}

In this section, a minimal working example of using \lstinline{ninfty} will be presented. We will, throughout this section and in \cref{sec:features}, assume that a finite group is our input. As already alluded to in \cref{sec:install}, all work should be done in the \lstinline{main.cpp} file, which is displayed below.

\begin{lstlisting}[language=C++]
#include "group_data/Cp2.h"
#include "ninfty.h"

int main() {
    dataSheet();
    return 0;
}
\end{lstlisting}

There are three lines of note here:
\begin{enumerate}[align=left]  
    \item[(Line 1)] In this line we select which group we wish to compute with. This points to a header file in the  \lstinline{group_data} folder. In this particular example, we are considering the group $C_{p^2}$ (specifically: $C_4$, but the choice of prime does not matter in this instance). This line can be changed to input any group, provided its data file has been generated and is in the \lstinline{group_data} folder.
    \item[(Line 2)] Here we are importing the \lstinline{ninfty.h} package to give us access to the algorithms provided by the library. This line should not be edited.
    \item[(Line 5)] This line calls the \lstinline{dataSheet()} function which takes no arguments. This will print out a collection of enumeration results regarding transfer systems for $C_{p^2}$.
\end{enumerate}

At this point, the code can be compiled and run, and this will differ from system to system. On a Unix machine with the \lstinline{g++} compiler installed one can run
\begin{lstlisting}[language=command.com,numbers=none]
g++ -O2 -pthread -std=c++20 main.cpp -o main.x
\end{lstlisting}
at the command line. Note the flag for selecting \lstinline{C++20}, the \lstinline{-pthread} flag which allows the code to run in parallel, and the \lstinline{-O2} flag for optimization. The compiled result is then run with
\begin{lstlisting}[language=command.com,numbers=none]
./main.x
\end{lstlisting}

This will output the data file for $C_{p^2}$ to the command line:
\begin{lstlisting}[language=HTML,numbers=none]
G=C4
#Transfer Systems=5
Complexity=2
Generation Statistics={1,3,1}
#Saturated Transfer Systems=4
Cosaturated Complexity=2
#Cosaturated Transfer Systems=4
Saturated Complexity=2
Width=2
#Flat transfers=4
#Premodel structures=13
#Composition closed structures=12
#Quillen structures=10
#Weak equivalence types=4
#Compatible pairs=12
\end{lstlisting}

If instead, we wanted to see this data sheet for the group $C_{pqr}$ (specifically: $C_{30}$), we need only edit line 1 of \lstinline{main.cpp} to
\begin{lstlisting}[language=C++]
#include "group_data/Cpqr.h"
\end{lstlisting}
and once again compile and run the executable and we obtain
\begin{lstlisting}[language=HTML,numbers=none]
G=C30
#Transfer Systems=450
Complexity=7
Generation Statistics={1,19,99,177,113,33,7,1}
#Cosaturated Transfer Systems=61
Cosaturated Complexity=4
#Saturated Transfer Systems=61
Saturated Complexity=4
Width=3
#Flat transfers=229
#Premodel structures=33903
#Composition closed structures=6949
#Quillen structures=1026
#Weak equivalence types=259
#Compatible pairs=13209
\end{lstlisting}

\section{Generating data files using Sage}\label{sec:Sage}

While many common groups are provided in the \lstinline{group_data} folder, it may be the case that a group which is not provided is required. In this case, one can use the included \lstinline{group_gen.sage} file. Of importance in this file are the first two lines:

\begin{lstlisting}[language=Python]
G0 = gap.SmallGroup(16,7)
group_name = "D8"
\end{lstlisting}

\begin{enumerate}[align=left]
    \item[(Line 1)] On this line we select the group that we are interested in. This has been designed so that GAP's small group library IDs can be used \cite{gap}.  For example,  \lstinline{SmallGroup(16,7)} corresponds to the dihedral group $D_8$. 
    \item[(Line 2)] This line allows for the naming of the output file. In this case the resulting file would be named ``\lstinline{D8.h}'' (note that the .h is not needed to be added by the user).
\end{enumerate}

To generate the corresponding data file, this Sage file need only be loaded in Sage in the usual fashion: by initiating Sage and navigating to the folder where the Sage file is located and running:
\begin{lstlisting}[language=command.com,numbers=none]
load("group_gen.sage")
\end{lstlisting}
This will produce some verbose output which is intended to allow the user to track the progress of the generation of the data. Once this is complete, the file will be created, and is immediately ready for use in  \lstinline{ninfty}.

Similarly, the \lstinline{lattice_gen.sage} file can be used to generate data files for arbitrary finite lattices. In this case we see the first two lines are similar to the \lstinline{group_gen.sage} file:

\begin{lstlisting}[language=Python]
L0 = Poset({0:[1,2,3,4,],1:[4,],2:[3,4,],3:[4,],4:[],})
lattice_name = "Pentagon"
\end{lstlisting}

The only difference here is that the variable \lstinline{L0} takes in a Sage poset, which here has been explicitly described using generators and relations in the usual format for Sage. It is assumed that the inputted poset is moreover a lattice.

\newpage

\section{Feature description}\label{sec:features}

We shall now outline all features of  \lstinline{ninfty}. As discussed in \cref{sec:context}, we have split these up into two categories \emph{basic} and \emph{advanced} depending on what the user wants to achieve with \lstinline{ninfty}. The advanced features are more granular in that they allow one to see if a given transfer system is saturated, for example, while the basic features would only tell the user how many saturated transfers there are.

\subsection{Basic features}\label{sec:featuresBASIC}

The basic features are simply functions that will print to the console a requested computation, in particular these functions take no inputs. In this section we will list these features, as well as reference where the concept is from. In places we shall also discuss how the algorithm is implemented into the code. Throughout, we will assume that we are working with a finite group $G$.

\hrulefill \vspace{3mm}

\begin{lstlisting}[language=C++,numbers=none]
printNumberOfTransfers()
\end{lstlisting}
\addcontentsline{toc}{subsubsection}{\lstinline{printNumberOfTransfers}}

This function will print the number of transfer systems for the group $G$ as defined in \cref{defn:xfer}, that is, it outputs $|\mathsf{Tr}(G)|$. As a test case we can run this for $C_{p^2}$ and we obtain 5, and for $C_{p^3}$ we obtain 14. Indeed, \cite{bbr} prove that the number of transfer systems for $C_{p^n}$ is given by the $(n+1)$-st Catalan number. Another test can that can be taken are rank two elementary abelian $p$-groups, where \cite{bao2023transfersystemsrankelementary} proves that there are exactly $2^{p+2}+p+1$ transfer systems for $G=(C_{p})^2$.

The algorithm that is used by \lstinline{ninfty} to do this computation is described in detail in \cite{abbswz}. In short, it uses Rubin's algorithm (\cite{rubin}) which finds the smallest transfer system containing a collection of pairs $(H,K)$. This algorithm is in fact a closure operator, and as such we use an efficient algorithm which finds all the closed sets of this closure operator, which are exactly the transfer systems. This algorithm has been implemented in parallel for efficiency.

The way that the algorithm is implemented means that we actually store all transfer systems in a variable \lstinline{ALL_STORE}, whose individual elements can be accessed. This variable is only computed once if required, and is used in many of the advanced features that will be discussed in \cref{sec:featuresADVANCED}.

\hrulefill \vspace{3mm}

\begin{lstlisting}[language=C++,numbers=none]
printNumberOfCosaturatedTransfers()
\end{lstlisting}
\addcontentsline{toc}{subsubsection}{\lstinline{printNumberOfCosaturatedTransfers}}

A transfer system $\sfT$ is \emph{cosaturated} (or disc-like) if it is generated (in the sense of Rubin's algorithm) by pairs $(H,G)$. These have been used, for example, in the study of splitting of incomplete rational Mackey functors \cite{barnes2024splittingrationalincompletemackey}. 

We use a similar algorithm as is used to generate all transfer systems, but with a modification that only allows the transfers to $G$ in the generation algorithm. Again, the collection of these is stored in the variable \lstinline{COSATURATED_STORE}. This function thus returns the size of this variable, that is, the number of cosaturated transfer systems. For $G=C_{pqr}$ we can check that there are 61 cosaturated transfer systems, while for $S_4$ there are 183.

\hrulefill \vspace{3mm}

\begin{lstlisting}[language=C++,numbers=none]
printNumberOfSaturatedTransfers()
\end{lstlisting}
\addcontentsline{toc}{subsubsection}{\lstinline{printNumberOfSaturatedTransfers}}

A \emph{saturated} transfer system is one that satisfies the 2-out-of-3 property. That is, for $L \leqslant K \leqslant H \leqslant G$, if two of $(L,K)$, $(L,H)$ and $(K,H)$ are in $\mathsf{T}$, then so is the third \cite{rubin}. This function returns the number of saturated transfer systems for $G$.

A study of this restricted class of transfer systems has been undertaken for some groups, for example in \cite{hmoo} and \cite{bao2023transfersystemsrankelementary}. While the cosaturated transfer systems correspond to disc-like operads, the saturated transfer systems share a close relation with the linear isometries operads for $G$. This connection is intricate, and has been studied in detail, for example in \cite{bannwart2023realizationsaturatedtransfersystems,arXiv:2311.08797}.

We have that for $G=C_{pqr}$ there are 61 saturated transfer systems. That is, there are as many saturated transfer systems for this group as there are cosaturated. This is not a coincidence, for abelian $G$ these numbers always coincide, and this is a consequences of a certain duality on the lattice of subgroups (see \cite{fooqw} and the discussion of the paragraph below). This is not true for non-abelian (specifically: non-Dedekind) groups. Indeed, there are only 132 saturated transfer systems for $S_4$.

The algorithm for computating this uses the observation that the saturated transfer systems for a lattice $L$ are in (non-canonical bijection) with the cosaturated transfer systems on $L^{\mathrm{op}}$ (\cite{bose2025combinatoricsfactorizationsystemslattices}). As such, we can simply run the generation algorithm used for the cosaturated transfer system on the opposite lattice of subgroups. Note, however, that by design that this does not actually produce the saturated transfer systems for $G$ itself, but just a (non-canonically) isomorphic set.

In the case that the collection of saturated transfer systems is required (as is the case for some of the advanced features) then a separate algorithm is used, \lstinline{saturatedTransfers} which checks the two-out-of-three property for all elements of \lstinline{ALL_STORE} and stores the output in \lstinline{SATURATED_STORE}. 

\hrulefill \vspace{3mm}

\begin{lstlisting}[language=C++,numbers=none]
printNumberOfUnderlyingTransfers()
\end{lstlisting}
\addcontentsline{toc}{subsubsection}{\lstinline{printNumberOfUnderlyingTransfers}}

When $\operatorname{Sub}(G)$ has no conjugation action (i.e., $G$ is a Dedekind group), then the transfer systems for $G$ coincide with the \emph{weak factorization systems} on the underlying lattice $\operatorname{Sub}(G)$.

However, in the case that $G$ has non-normal subgroups, then the conjugation condition in the definition of a transfer system means that it is far more restrictive to be a transfer system for $G$ than a weak factorization system on the underlying lattice of $\operatorname{Sub}(G)$. This function will return the number of weak factorization systems on this underlying lattice. The relation between these two counts has been considered (albeit from a slightly different point of view) in \cite{BMO_lift}. We can check, for example, that while $S_3$ has $9$ transfer systems, that there are $36$ weak factorization systems on the underlying non-equivariant lattice of subgroups.

To facilitate this, a further modified version of Rubin's algorithm is used which skips the step where the conjugation is added. The underlying transfer systems are stored in the variable \lstinline{UNDERLYING_STORE}.

\hrulefill \vspace{3mm}

\begin{lstlisting}[language=C++,numbers=none]
printNumberOfConjugacyTransfers()
\end{lstlisting}
\addcontentsline{toc}{subsubsection}{\lstinline{printNumberOfConjugacyTransfers}}

Following on from the previous function, we could instead consider the number of transfer systems on $\mathrm{Sub}(G)/G$. While $\mathrm{Sub}(G)/G$ need not be a lattice in general, we can use a more general definition of a transfer system which works for posets \cite{BMO_lift}. This function returns the number of transfer systems on this poset, and these are stored in the variable \lstinline{CONJUGACY_STORE}.

Again, while $S_3$ has 9 transfer systems, we have that $\mathrm{Sub}(S_3)/S_3 \cong [1] \times [1]$ and this has 10 transfer systems. The theory of \emph{liftable} transfer systems from \cite{BMO_lift} describes how we can obtain the transfer systems for $S_3$ from these 10. Note that this comparison can only be made when $G$ is a \emph{lossless} group.

\hrulefill \vspace{3mm}

\begin{lstlisting}[language=C++,numbers=none]
printWidth()
\end{lstlisting}
\addcontentsline{toc}{subsubsection}{\lstinline{printWidth}}

While Rubin's algorithm allows us to build the smallest transfer system containing a given collection of pairs $(H,K)$, we can consider running this algorithm in reverse, which allows one to obtain a \emph{minimal basis} for any given transfer system $\mathsf{T}$. Note that we disregard any trivial pairs $(H,H)$. The sizes of these bases form the main thrust of investigation in \cite{abbswz}.

We define the \emph{width} of a finite group $G$ to be the number of non-trivial generators in a minimal generating set for the complete transfer systems (that is, the transfer system with all possible pairs $(H,K)$). It was  proved in \cite{abbswz} that the width of $G$ coincides with the number of meet-irreducible subgroups of $G$. These are those proper subgroups $H < G$ such that $H$ cannot be written as the intersection of two other proper subgroups. In particular, every maximal subgroup is meet-irreducible, but the converse may not hold. For example, for $C_{p^2q}$ the width is 3 while there are only two maximal subgroups.

To compute the width effectively, we can run Rubin's algorithm in reverse on the complete transfer system. This means that we do not need to compute all transfer systems to be able to compute the width.

\hrulefill \vspace{3mm}

\begin{lstlisting}[language=C++,numbers=none]
printComplexity()
\end{lstlisting}
\addcontentsline{toc}{subsubsection}{\lstinline{printComplexity}}

While the width computes the size of a minimal basis for the complete transfer system, the \emph{complexity} of a group is the largest sized basis over all possible transfer systems \cite{abbswz}. What may be surprising at first pass is that this number is, in general, much larger than the width of $G$. Indeed, it is proved in \cite{abbswz} that the complexity of $C_{pqr}$ is 7, while the width is 3.

Unlike the computation of the width, it is not possible \textit{a priori} to compute the complexity without computing all transfer systems. However, this number is naturally computed when all transfer systems are computed as it is exactly the number of iterations that he implemented algorithm needs to run for.

\hrulefill \vspace{3mm}

\begin{lstlisting}[language=C++,numbers=none]
printNumberOfMaximallyGenerated()
\end{lstlisting}
\addcontentsline{toc}{subsubsection}{\lstinline{printNumberOfMaximallyGenerated}}

Following on from the computation of the complexity, this function prints out the number of transfer systems which have minimal basis size the complexity. This seems to be an interesting question which has not yet seen much investigation. Surprisingly many groups have only one transfer system which realises the complexity. For example, while $S_4$ has 14 maximally generated transfer systems, $S_5$ only has 1.

These transfer systems are stored in the variable \lstinline{MAXIMALLY_GENERATED}.

\hrulefill \vspace{3mm}

\begin{lstlisting}[language=C++,numbers=none]
printNumberOfFlatTransfers()
\end{lstlisting}
\addcontentsline{toc}{subsubsection}{\lstinline{printNumberOfFlatTransfers}}

Let $\mathsf{T}$ be a transfer system for $G$. Then there is a unique normal subgroup $F$ which is minimal among all subgroups $H$ such that the pair $(H,G)$ is in $\sfT$. This $F$ is referred to as the \emph{minimal fibrant subgroup} in \cite{BMO_enumeration} due to its relation with model structures (to be discussed shortly).

A transfer system is defined to be \emph{flat} if $\mathrm{T}$ restricted to $\mathsf{Sub}(F)$ is trivial. That is, there are no non-trivial pairs $(H,K)$ for $H \leqslant K \leqslant F$. These were first studied in \cite{MR4565666} where it was shown that free incomplete Tambara functors are almost never flat. For the family of cyclic groups $C_{p^n}$ the number of flat transfer systems are counted by partial sums of Catalan numbers, and in particular, as $n \to \infty$ the number of flat transfer systems approaches $1/3$ of all transfer systems.

This function will return the number of flat transfer systems for $G$ and moreover stores the flat transfer systems in the variable \lstinline{FLAT_STORE}.

\newpage
\hrulefill \vspace{3mm}

\begin{lstlisting}[language=C++,numbers=none]
printNumberOfTransferPairs()
\end{lstlisting}
\addcontentsline{toc}{subsubsection}{\lstinline{printNumberOfTransferPairs}}

We now move into the realm of looking at certain pairs of transfer systems. It was proved in \cite{fooqw} that $\mathsf{Tr}(G)$ is a complete lattice. That is, we can consider meets and joins of transfer systems. The goal of this function is to compute the number of intervals in this lattice, that is, the number of (potentially trivial) inclusions $\sfT \leqslant \sfT'$ in $\mathsf{Tr}(G)$. These intervals are stored in \lstinline{TRANSFER_LATTICE} whose objects are pairs of indices $(i,j)$ such that \lstinline{ALL_STORE[i]} $\leqslant$ \lstinline{ALL_STORE[j]}. 

These intervals correspond to \emph{premodel structures} on the lattice $\mathrm{Sub}(G)$ in the sense of \cite{barton, BOOR}. In more detail, any complete lattice is in particular a complete category, and as such we can talk about \emph{weak factorization systems} on this lattice. It was proved in \cite{fooqw} that for Dedekind groups that there is a bijection between transfer systems for $G$ and weak factorization systems on the lattice $\operatorname{Sub}(G)$. The edges which appear in the transfer system correlate to the \textit{fibrations} of this weak factorization system. Assocaited to this are those morphisms which have the left lifting property with respect to the fibrations, and this gives us a class of \emph{acyclic cofibrations}.

In the case that we have an interval of transfer systems $\sfT \leqslant \sfT'$ we can take $\sfT$ to be the fibrations and $\sfT'$ to be the \emph{acyclic fibrations} of a \emph{premodel structure}. This allows us to define a notion of \emph{weak equivalence} to be those morphisms which can be written as the composition of an acyclic cofibration and an acyclic fibration. 

This function returns the number of premodel structures on the subgroup lattice of $G$. Note that in the case that $G$ is non-Dedekind then this will return the number of premodel structures on the underlying lattice of $\mathrm{Sub}(G)$ which can be assembled out of the transfer systems for $G$. Work of \cite{bbr} tells us that the number of premodel structures in the case of $C_{p^n}$ are given exactly by the number of intervals in the $(n+1)$-Tamari lattice, which has been computed by \cite{chapoton}.

\hrulefill \vspace{3mm}

\begin{lstlisting}[language=C++,numbers=none]
printNumberOfCClosedPairs()
\end{lstlisting}
\addcontentsline{toc}{subsubsection}{\lstinline{printNumberOfCClosedPairs}}

In \cite{BMO_comp}, the notion of \emph{composition closed} premodel structures was introduced. These are those premodel structures such that the collection of weak equivalences as defined above is closed under composition. This function returns the number of composition closed premodel structures on $\operatorname{Sub}(G)$, with the same caveat about the non-Dedekind setting as above.

For $G = C_{p^n}$ it is possible to identify these composition closed model structures. Instead of intervals of the Tamari lattice, we  need only consider intervals of a coarsening of the Tamari lattice, the Kreweras lattice \cite {BMO_comp}.

\hrulefill \vspace{3mm}

\begin{lstlisting}[language=C++,numbers=none]
printNumberOfQuillenPairs()
\end{lstlisting}
\addcontentsline{toc}{subsubsection}{\lstinline{printNumberOfQuillenPairs}}

Among the composition closed premodel structures there are those such that the weak equivalences satisfy the stronger property of being closed under the two-out-of-three property. These are exactly the \emph{Quillen model structures} on the subgroup lattice in the sense of \cite{quillen, handbook}, and it is the number of these that this function returns.

In \cite{BOOR} the number of Quillen model structures on $C_{p^n}$ was computed as $\binom{2n+1}{n}$. Interestingly the methods used in this paper seem not to extend beyond this case.

\newpage

\hrulefill \vspace{3mm}

\begin{lstlisting}[language=C++,numbers=none]
printNumberOfWeakEquivalenceTypes()
\end{lstlisting}
\addcontentsline{toc}{subsubsection}{\lstinline{printNumberOfWeakEquivalenceTypes}}

Given a model structure on $\mathrm{Sub}(G)$, it is (usually) the case that there are other model structures which have the same weak equivalence type. That is, the weak equivalences are the same but the fibrations differ. This function returns the number of unique weak equivalence types among all model structures. This is helpful when determining the amount of possible homotopy types that can be obtained from the lattice in question.

\hrulefill \vspace{3mm}

\begin{lstlisting}[language=C++,numbers=none]
printNumberOfCompatiblePairs()
\end{lstlisting}
\addcontentsline{toc}{subsubsection}{\lstinline{printNumberOfCompatiblePairs}}

 Algebraically, transfer systems control \emph{incomplete Tambara functors} \cite{MR3773736}. If we have two transfer systems $\mathsf{T}$ and $\mathsf{T}'$ such that $\mathsf{T} \leqslant \mathsf{T}'$ then we say that these transfer systems are \emph{compatible} if whenever $A \leqslant G$ and $B,C \leqslant A$ are subgroups such that $(B,A) \in \mathsf{T}$ and $((B \cap C) ,B) \in \mathsf{T}'$ then $(C, A) \in \mathsf{T}'$ \cite{Chan_tambara}. These compatible pairs control \emph{bi-incomplete Tambara functors} \cite{BH_bitambara}. This function returns the number of bi-incomplete Tambara functors for $G$ by computing the number of compatible pairs.

Compatible pairs have been studied in \cite{hill2022countingcompatibleindexingsystems} and \cite{mazur2024uniquelycompatibletransfersystems}, with the former proving that the number of compatible pairs for $C_{p^n}$ is counted by the Fuss--Catalan numbers $A_{n+1}(3,1)$, which (unique to $C_{p^n}$) coincides with the number of composition closed model structures. Note that it is not necessarily the case that a given transfer systems is compatible with itself.

\hrulefill \vspace{3mm}

\begin{lstlisting}[language=C++,numbers=none]
printSageTransferPoset()
\end{lstlisting}
\addcontentsline{toc}{subsubsection}{\lstinline{printSageTransferPoset}}

As mentioned previously, it was proved in \cite{fooqw} that the collection of all transfer systems is itself a lattice. This function will return a Sage command for this poset so that it can be investigated further.

\hrulefill \vspace{3mm}

\begin{lstlisting}[language=C++,numbers=none]
printSageCClosedPoset()
\end{lstlisting}
\addcontentsline{toc}{subsubsection}{\lstinline{printSageCClosedPoset}}

Similar to the above, this function will return a Sage string for the poset of transfer systems where $\mathsf{T} \leqslant \sfT'$ if and only if this pair forms a composition closed model structure.

\hrulefill \vspace{3mm}

\begin{lstlisting}[language=C++,numbers=none]
printSageQuillenPoset()
\end{lstlisting}
\addcontentsline{toc}{subsubsection}{\lstinline{printSageQuillenPoset}}

Similar to the above, this function will return a Sage string for the poset of transfer systems where $\mathsf{T} \leqslant \sfT'$ if and only if this pair forms a Quillen model structure.

\newpage

\hrulefill \vspace{3mm}

\begin{lstlisting}[language=C++,numbers=none]
dataSheet()
dataSheetLatex()
\end{lstlisting}
\addcontentsline{toc}{subsubsection}{\lstinline{dataSheet} and \lstinline{dataSheetLatex}}

As already discussed in \cref{sec:MWE}, the function \lstinline{dataSheet} will print out a summary of all the enumerative output from the functions above. The function \lstinline{dataSheetLatex} is identical in form to this, but instead returns the LaTeX code for this information to be formatted in a table, and has been designed for easy inclusion into papers. 
For example, for $G=A_5$ the code yields \cref{table:a5}.
\begin{table}[h]
\begin{tabular}{|cc|}
\hline
\multicolumn{2}{|c|}{$G = A_5$} \\ \hline
\multicolumn{1}{|c|}{\#Transfer systems} & 987\\ \hline
\multicolumn{1}{|c|}{Complexity} & 8\\ \hline
\multicolumn{1}{|c|}{Width} & 5\\ \hline
\multicolumn{1}{|c|}{Generation values} & \{1,23,126,285,308,175,57,11,1\}\\ \hline
\multicolumn{1}{|c|}{\#Saturated} & 55\\ \hline
\multicolumn{1}{|c|}{Saturated complexity} & 5\\ \hline
\multicolumn{1}{|c|}{\#Cosaturated} & 61\\ \hline
\multicolumn{1}{|c|}{Cosaturated complexity} & 5\\ \hline
\multicolumn{1}{|c|}{\#Flat} & 450\\ \hline
\multicolumn{1}{|c|}{\#Premodel structures} & 151816\\ \hline
\multicolumn{1}{|c|}{\#C.closed structures} & 25874\\ \hline
\multicolumn{1}{|c|}{\#Quillen structures} & 1813\\ \hline
\multicolumn{1}{|c|}{\#Weak equivalence types} & 445\\ \hline
\multicolumn{1}{|c|}{\#Compatible pairs} & 49651 \\ \hline
\end{tabular}
\caption{Data sheet for $A_5$.}\label{table:a5}
\end{table}

\hrulefill \vspace{3mm}
\begin{lstlisting}[language=C++,numbers=none]
printSubgroupDictionary()
\end{lstlisting}
\addcontentsline{toc}{subsubsection}{\lstinline{printSubgroupDictionary}}

This is a utility function which will print out the subgroup dictionary of $G$ as stored by \lstinline{ninfty}. In the case that $G$ is not a Dedekind group it will also print out the conjugacy classes of the subgroups.  For example, in the case that $G=S_3$ we get:
\begin{lstlisting}[language=HTML,numbers=none]
{0:1}
{1:C2(1)}
{2:C2(2)}
{3:C2(3)}
{4:C3}
{5:S3}

Conjugacy Classes:
[0]
[1,2,3]
[4]
[5]
\end{lstlisting}
Unwrapping this, we see that $S_3$ has six subgroups, the 0th subgroup is the trivial subgroup, subgroups 1,2, and 3 are all copies of $C_2$, the fourth subgroup is a $C_3$ and then we have the whole group. The conjugacy class data moreover informs us that subgroups 1,2, and 3 are all in one conjugacy class, while subgroup 4 (the copy of $C_3$) is a normal subgroup as we would expect.

\hrulefill \vspace{3mm}

\begin{lstlisting}[language=C++,numbers=none]
printAllTransfers()
\end{lstlisting}
\addcontentsline{toc}{subsubsection}{\lstinline{printAllTransfers}}

This function will print out all the transfer systems for $G$. Running this for $S_3$ we get:
\begin{lstlisting}[language=HTML,mathescape=true]
{$\varnothing$}
{(0,4)}
{(0,1),(0,2),(0,3)}
{(0,1),(0,2),(0,3),(4,5)}
{(0,1),(0,2),(0,4),(0,3)}
{(0,1),(0,2),(0,4),(0,3),(0,5)}
{(0,1),(0,2),(0,4),(0,3),(0,5),(4,5)}
{(0,1),(0,2),(0,4),(0,3),(0,5),(1,5),(2,5),(3,5)}
{(0,1),(0,2),(0,4),(0,3),(0,5),(1,5),(2,5),(4,5),(3,5)}
\end{lstlisting}
We see that there are 9 transfer systems as expected. The first transfer system is the trivial one. The second transfer system has only one non-trivial pair $(0,4)$, which, referring back to the subgroup dictionary is the relation $(e,C_3)$.

\hrulefill \vspace{3mm}

\subsection{Advanced features}\label{sec:featuresADVANCED}

We now move onto the discussion of (some of) the more advanced features that are provided in \lstinline{ninfty}. Many of these are designed so that the user can probe features of individual transfer systems.  Already in \cref{sec:featuresBASIC} we have discussed how calling some of the basic features populates certain variables which these algorithms can then be run on. The discussions of these algorithms will be briefer than what was discussed in \cref{sec:featuresBASIC} as the theory has already been covered there. Of course, by combining these algorithms far more bespoke calculations can be performed.

To understand these algorithms, it is worth discussing how a collection of pairs $(H,K)$ is stored in \lstinline{ninfty}. In particular, this is how transfer systems are stored.

The first variable to understand is \lstinline{lattice} in the group data file. This is a vector consisting of pairs of unsigned integers $(i,j)$. This corresponds to the inclusion \lstinline{subgroup_dictionary[i]} to \lstinline{subgroup_dictionary[j]}. \lstinline{ninfty} then stores a collection of pairs $(H,K)$ as a \lstinline{std::vector<unsigned>}, that is, a vector of unsigned integers. Each $i$ in this structure corresponds to the $i$th entry of the \lstinline{lattice} variable.

\hrulefill \vspace{3mm}

\begin{lstlisting}[language=C++,numbers=none]
std::vector<std::vector<unsigned>> transferFind(bool verbose, enum gen_type)
\end{lstlisting}
\addcontentsline{toc}{subsubsection}{\lstinline{transferFind}}

This function is the main function which is used to populate the various variables of transfer systems. It takes two optional arguments:
\begin{enumerate}
\item The first is a boolean \lstinline{verbose}, which is default set to \lstinline{false}, is a toggle on if the generation data is outputted or not while running the algorithm. Having this set to \lstinline{true} can help to track the progress of the algorithm.
\item The second argument is an \lstinline{enum} which picks the generation type. The possible options for this variable for this are 

\centerline{\lstinline{ALL, SATURATED, COSATURATED, UNDERLYING, CONJUGACY}}

with the default being \lstinline{ALL}. These correspond to the different types of generation that have already been covered in \cref{sec:featuresBASIC}.
\end{enumerate} 
The output of this function is a vector which contains all of the requested transfer systems, which are moreover stored in the variables already discussed in \cref{sec:featuresBASIC}.

\newpage

\hrulefill \vspace{3mm}

\begin{lstlisting}[language=C++,numbers=none]
bool isSaturated(std::vector<unsigned> rhs)
\end{lstlisting}
\addcontentsline{toc}{subsubsection}{\lstinline{isSaturated}}

This function returns \lstinline{true} if the inputted transfer system \lstinline{rhs} is saturated, and \lstinline{false} if not.

\hrulefill \vspace{3mm}

\begin{lstlisting}[language=C++,numbers=none]
bool isCosaturated(std::vector<unsigned> rhs)
\end{lstlisting}
\addcontentsline{toc}{subsubsection}{\lstinline{isCosaturated}}

This function returns \lstinline{true} if the inputted transfer system \lstinline{rhs} is cosaturated, and \lstinline{false} if not.

\hrulefill \vspace{3mm}

\begin{lstlisting}[language=C++,numbers=none]
bool isFlat(std::vector<unsigned> rhs)
\end{lstlisting}
\addcontentsline{toc}{subsubsection}{\lstinline{isFlat}}

This function returns \lstinline{true} if the inputted transfer system \lstinline{rhs} is flat, and \lstinline{false} if not.

\hrulefill \vspace{3mm}

\begin{lstlisting}[language=C++,numbers=none]
std::vector<unsigned> leftSet(std::vector<unsigned> rhs)
\end{lstlisting}
\addcontentsline{toc}{subsubsection}{\lstinline{leftSet}}

This function returns the left set (in the sense of weak factorization systems) of the inputted transfer system \lstinline{rhs} using \cite[Proposition 4.15]{fooqw}. This is again stored in a \lstinline{std::vector<unsigned>} in the same aforementioned format.

\hrulefill \vspace{3mm}

\begin{lstlisting}[language=C++,numbers=none]
std::vector<unsigned> saturatedHull(std::vector<unsigned> rhs)
\end{lstlisting}
\addcontentsline{toc}{subsubsection}{\lstinline{saturatedHull}}

Every transfer system has a unique (minimal among all saturated transfer systems) saturated hull in which it is contained in. One can obtain this transfer system by simply closing up under the 2-out-of-3 property. This function returns the saturated hull of the inputted transfer system \lstinline{rhs}.

\hrulefill \vspace{3mm}

\begin{lstlisting}[language=C++,numbers=none]
std::vector<unsigned> cosaturatedCore(std::vector<unsigned> rhs)
\end{lstlisting}
\addcontentsline{toc}{subsubsection}{\lstinline{cosaturatedCore}}

Every transfer system has a unique (maximal among all cosaturated transfer systems) cosaturated hull which is contained in it. One can obtain this transfer system by simply completing the collection of pairs $(H,G)$ appearing in \lstinline{rhs} to a transfer system.  This function returns the cosaturated core of the inputted transfer system \lstinline{rhs}.

\hrulefill \vspace{3mm}

\begin{lstlisting}[language=C++,numbers=none]
std::vector<unsigned> dualTransferSystem(std::vector<unsigned> rhs)
\end{lstlisting}
\addcontentsline{toc}{subsubsection}{\lstinline{dualTransferSystem}}

In the case that $G$ is a cyclic group, then it was proved in \cite{fooqw} that $\mathsf{Tr}(G)$ admits a self-duality. This crucially uses the canonical self-duality on the subgroup lattice of a cyclic group.

This function will return the dual of the transfer system \lstinline{rhs} under the caveat that the group in question is a cyclic group. When $G$ is not cyclic, this function will simply return \lstinline{rhs}.

\newpage

\hrulefill \vspace{3mm}

\begin{lstlisting}[language=C++,numbers=none]
unsigned minimalFibrantSubgroup(std::vector<unsigned> rhs)
\end{lstlisting}
\addcontentsline{toc}{subsubsection}{\lstinline{minimalFibrantSubgroup}}

We have already discussed how there is a unique minimal subgroup $H$ such that $(H,G)$ is in the tranfer system, this is the \emph{minimal fibrant subgroup} of the transfer system. This function returns the index in \lstinline{subgroup_dictionary} of the minimal fibrant subgroup of the inputted transfer system \lstinline{rhs}.

\hrulefill \vspace{3mm}

\begin{lstlisting}[language=C++,numbers=none]
std::vector<unsigned> findBasis(std::vector<unsigned> rhs)
\end{lstlisting}
\addcontentsline{toc}{subsubsection}{\lstinline{findBasis}}

This function returns a minimal basis for the inputted transfer system \lstinline{rhs}. This is a set of pairs $(H,K)$ which will generate \lstinline{rhs} under Rubin's algorithm with no redundancies. Note that this minimal basis is not unique, but the code uses a deterministic algorithm (namely: running Rubin's algorithm in reverse) and thus will always return the same generating set. Moreover, any minimal generating set will have the same cardinality as this one. 

\hrulefill \vspace{3mm}

\begin{lstlisting}[language=C++,numbers=none]
std::vector<unsigned> weakEquivalences(std::vector<unsigned> AF,std::vector<unsigned> F)
\end{lstlisting}
\addcontentsline{toc}{subsubsection}{\lstinline{weakEquivalences}}

This function is crucial in the functions for finding the various types of model structures arising from the transfer systems. In particular, it returns the weak equivalences of the given input data. The input data is two collections of pairs, namely \lstinline{AF}, the \emph{acyclic fibrations}, and \lstinline{F}, the \emph{fibrations}. We are implicitly assuming that \lstinline{AF} $\subseteq$ \lstinline{F}.

\hrulefill \vspace{3mm}

\begin{lstlisting}[language=C++,numbers=none]
unsigned modelCheck(std::vector<unsigned> AF,std::vector<unsigned> F)
\end{lstlisting}
\addcontentsline{toc}{subsubsection}{\lstinline{modelCheck}}

The inputs for this function are the same as the inputs for \lstinline{weakEquivalences}, namely a set of acyclic fibrations and a set of fibrations. The return of this function has three options:
\begin{itemize}
    \item 2 if the pair is a Quillen model structure,
    \item 1 if the pair is a composition closed model structure which is not Quillen,
    \item 0 if the pair is not a composition closed model structure.
\end{itemize}

\hrulefill \vspace{3mm}

\begin{lstlisting}[language=C++,numbers=none]
bool isCompatible(std::vector<unsigned> transfer_m, std::vector<unsigned> transfer_a)\end{lstlisting}
\addcontentsline{toc}{subsubsection}{\lstinline{isCompatible}}

This function will check if two transfer systems \lstinline{transfer_m} $\leqslant$ \lstinline{transfer_a} form a compatible pair or not.

\hrulefill \vspace{3mm}

\begin{lstlisting}[language=C++,numbers=none]
edgesToTikz(std::vector<unsigned> rhs)
\end{lstlisting}
\addcontentsline{toc}{subsubsection}{\lstinline{edgesToTikz}}

This is a utility function which will print out a TikZ string to the console which displays the collection of pairs $(H,K)$ in \lstinline{rhs} (which need not be a transfer system!).

At the base level, this will plot the edges of the transfer system on the poset $\operatorname{Sub}(G)/G$ (which may not be a perfect representation in a lossy group, but it is unfeasible to plot all edges on the entire subgroup lattice). The conjugacy classes of subgroups  will be laid out in a circle with labels on the nodes given by the entries of \lstinline{subgroup_dictionary}.

Let us look at an example. The following is a particular transfer system on $Q_8$:
\begin{figure}[h]
\centering
\begin{tikzpicture}[scale=0.5]
\node[inner sep=0cm] (0) at (0.0000,3.0000) {$1$};
\node[inner sep=0cm] (1) at (2.5980,1.5000) {$C2$};
\node[inner sep=0cm] (2) at (2.5980,-1.499) {$C4$};
\node[inner sep=0cm] (3) at (0.0000,-3.000) {$C4$};
\node[inner sep=0cm] (4) at (-2.598,-1.500) {$C4$};
\node[inner sep=0cm] (5) at (-2.598,1.4999) {$Q8$};
\draw[black!10,->] (0) edge (1);
\draw[black!10,->] (0) edge (2);
\draw[black!10,->] (0) edge (3);
\draw[black!10,->] (0) edge (4);
\draw[black!10,->] (0) edge (5);
\draw[black!10,->] (3) edge (5);
\draw[codepurple,->] (1) edge (2);
\draw[codepurple,->] (1) edge (3);
\draw[codepurple,->] (1) edge (4);
\draw[codepurple,->] (1) edge (5);
\draw[codepurple,->] (2) edge (5);
\draw[codepurple,->] (4) edge (5);
\end{tikzpicture}
\end{figure}

While this figure does give a representation of a transfer system, it can be improved somewhat. The group data file for each group has optional variables which can be manually edited to make these TikZ diagrams more aesthetic, as we will describe now.

First, in the group data file is the variable \lstinline{pretty_subgroup_dictionary} which is a vector of strings. This can be used to relabel the vertices. The code checks that a name has been given to each conjugacy class of subgroup, this is done in the order in the above diagram, starting at $1$ and going around the subgroups clockwise. For $Q_8$ we can use:

\begin{lstlisting}[language=C++,numbers=none]
std::vector<std::string> pretty_subgroup_dictionary{
    "1",
    "C_2",
    "C_4",
    "C_4",
    "C_4",
    "Q_8"
};
\end{lstlisting}

Rerunning the \lstinline{edgesToTikz} code with this variable updated then produces the following diagram:

\begin{figure}[h]
\centering
\begin{tikzpicture}[scale=0.5]
\node[inner sep=0cm] (0) at (0.0000,3.0000) {$1$};
\node[inner sep=0cm] (1) at (2.5980,1.5000) {$C_2$};
\node[inner sep=0cm] (2) at (2.5980,-1.499) {$C_4$};
\node[inner sep=0cm] (3) at (0.0000,-3.000) {$C_4$};
\node[inner sep=0cm] (4) at (-2.598,-1.500) {$C_4$};
\node[inner sep=0cm] (5) at (-2.598,1.4999) {$Q_8$};
\draw[black!10,->] (0) edge (1);
\draw[black!10,->] (0) edge (2);
\draw[black!10,->] (0) edge (3);
\draw[black!10,->] (0) edge (4);
\draw[black!10,->] (0) edge (5);
\draw[black!10,->] (3) edge (5);
\draw[codepurple,->] (1) edge (2);
\draw[codepurple,->] (1) edge (3);
\draw[codepurple,->] (1) edge (4);
\draw[codepurple,->] (1) edge (5);
\draw[codepurple,->] (2) edge (5);
\draw[codepurple,->] (4) edge (5);
\end{tikzpicture}
\end{figure}

Next we change the positioning of the nodes. This can be adjusted using the \lstinline{vertex_layout} variable in the group data file. In the same ordering as the \lstinline{pretty_subgroup_dictionary} we can provide coordinates for each conjugacy class in the diagram:

\begin{lstlisting}[language=C++,numbers=none]
std::vector<std::string> vertex_layout{
    "(2,0)",
    "(2,0.803)",
    "(2,1.76)",
    "(3.88,1.76)",
    "(0.125,1.76)",
    "(2,2.71)"
};
\end{lstlisting}

\newpage
This now results in the following updated TikZ diagram:

\begin{figure}[!h]
\centering
\begin{tikzpicture}
\node[inner sep=0cm] (0) at (2,0){$1$};
\node[inner sep=0cm] (1) at (2,0.803){$C_2$};
\node[inner sep=0cm] (2) at (2,1.76){$C_4$};
\node[inner sep=0cm] (3) at (3.88,1.76){$C_4$};
\node[inner sep=0cm] (4) at (0.125,1.76){$C_4$};
\node[inner sep=0cm] (5) at (2,2.71){$Q_8$};
\draw[black!10,->] (0) edge (1);
\draw[black!10,->] (0) edge (2);
\draw[black!10,->] (0) edge (3);
\draw[black!10,->] (0) edge (4);
\draw[black!10,->] (0) edge (5);
\draw[black!10,->] (3) edge (5);
\draw[codepurple,->] (1) edge (2);
\draw[codepurple,->] (1) edge (3);
\draw[codepurple,->] (1) edge (4);
\draw[codepurple,->] (1) edge (5);
\draw[codepurple,->] (2) edge (5);
\draw[codepurple,->] (4) edge (5);
\end{tikzpicture}
\end{figure}

This updated diagram now has the issue that the edge, for example, from $C_2$ to $Q_8$ covers up the edges beneath it. This issue can be avoided by bending the edges to avoid these overlaps. This can be achieved using the final optional variable, \lstinline{edge_options}, which takes the form of an upper-triangular $n \times n$ vector where entry $[i][j]$ contains edge options for the edge from the $i$th subgroup to the $j$th subgroup. In this case we use the following options:

\begin{lstlisting}[language=C++,numbers=none]
std::vector<std::vector<std::string>> edge_options{
    {"","","[bend right]","","","[bend left]"},
    {"","","","","","[bend left]"},
    {"","","","","",""},
    {"","","","","",""},
    {"","","","","",""},
    {"","","","","",""},
};
\end{lstlisting}

\begin{figure}
\centering
\includegraphics[scale=0.51]{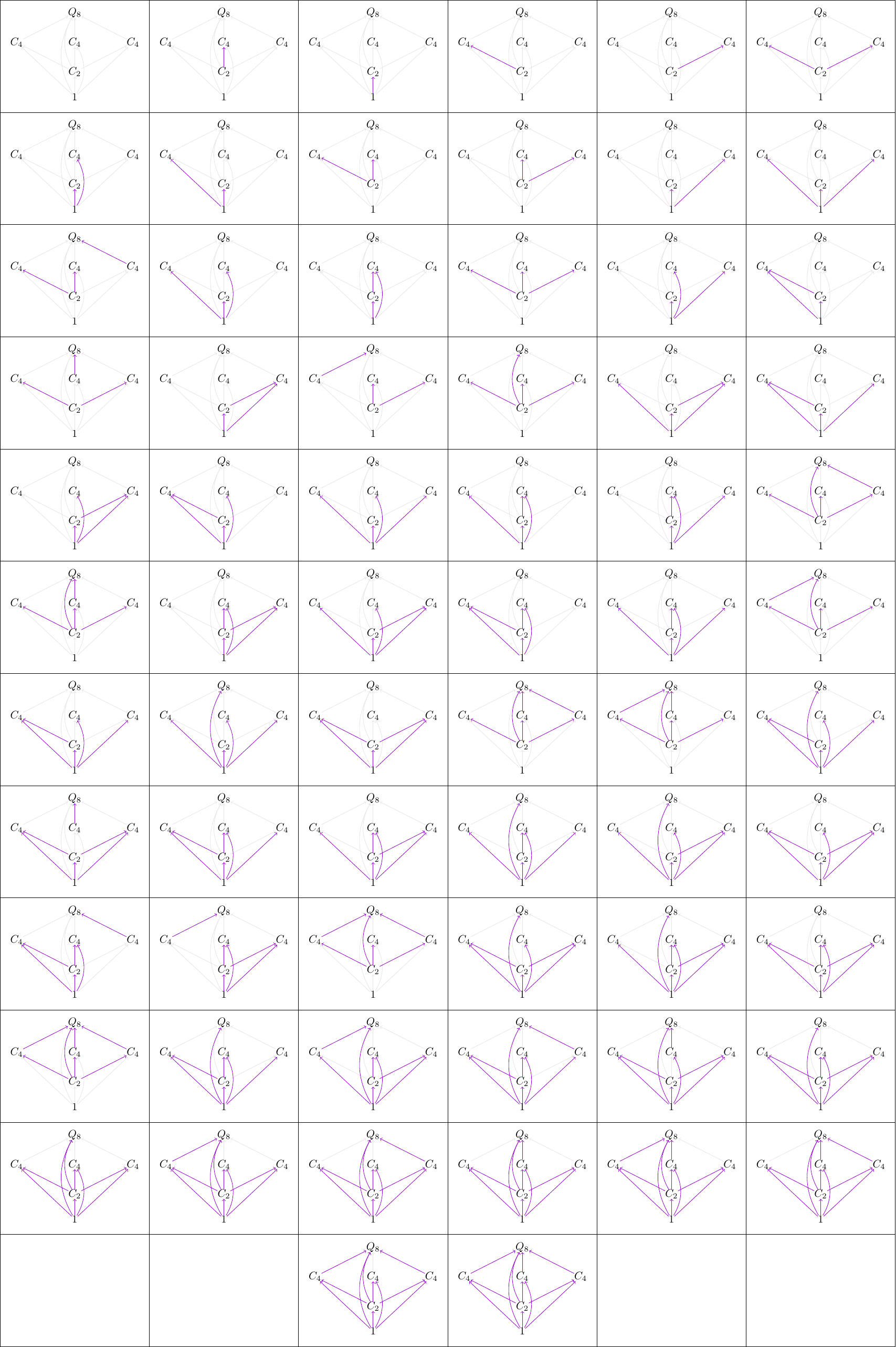}
\vspace*{-5mm}
\caption{The 68 transfer systems for $Q_8$.}\label{fig:q8}
\end{figure}

This results in our final layout and setup for $Q_8$:

\begin{figure}[h]
\centering
\begin{tikzpicture}
\node[inner sep=0cm] (0) at (2,0){$1$};
\node[inner sep=0cm] (1) at (2,0.803){$C_2$};
\node[inner sep=0cm] (2) at (2,1.76){$C_4$};
\node[inner sep=0cm] (3) at (3.88,1.76){$C_4$};
\node[inner sep=0cm] (4) at (0.125,1.76){$C_4$};
\node[inner sep=0cm] (5) at (2,2.71){$Q_8$};
\draw[black!10,->] (0) edge (1);
\draw[black!10,->] (0) edge[bend right] (2);
\draw[black!10,->] (0) edge (3);
\draw[black!10,->] (0) edge (4);
\draw[black!10,->] (0) edge[bend left] (5);
\draw[black!10,->] (3) edge (5);
\draw[codepurple,->] (1) edge (2);
\draw[codepurple,->] (1) edge (3);
\draw[codepurple,->] (1) edge (4);
\draw[codepurple,->] (1) edge[bend left] (5);
\draw[codepurple,->] (2) edge (5);
\draw[codepurple,->] (4) edge (5);
\end{tikzpicture}
\end{figure}
After setting this up, it becomes a near triviality to print out and display all transfer systems for $Q_8$, as done in \cref{fig:q8}.

\section{Support}

Any questions or suggestions should be directed to \url{s.balchin@qub.ac.uk}. The \lstinline{GitHub} repository should always be checked for the most up to date version of \lstinline{ninfty} and this documentation.

\newpage

\newcommand{\etalchar}[1]{$^{#1}$}

\end{document}